\documentclass[reqno]{amsart}
\usepackage[margin=3.5cm]{geometry}
\usepackage[dvipsnames]{xcolor}
\definecolor{darkblue}{rgb}{0,0,0.7}
\definecolor{darkred}{rgb}{0.7,0,0}
\usepackage[T1]{fontenc}
\usepackage{amsfonts, amsthm, amssymb,verbatim,amscd,amsmath, blindtext}

\usepackage{multirow}
\usepackage{graphicx}
\usepackage{enumitem,yhmath}
\usepackage{amssymb}
\usepackage{ytableau}
\usepackage{tikz}
\usepackage{tikz-cd}
\usepackage{float, pifont}
\usepackage{mathtools}
\usepackage[colorlinks=true, linkcolor=darkred, citecolor=darkblue, urlcolor=blue, pagebackref=true, breaklinks=true]{hyperref}
\usepackage[capitalise]{cleveref}

\DeclareMathOperator{\lcm}{lcm}
\DeclareMathOperator{\mingens}{Mingens}

\DeclareMathOperator{\pdim}{pdim}

\DeclareMathOperator{\HS}{HS}

\newcommand{\TT}{\mathbb{T}}

\newtheorem{theorem}{Theorem}[section]
\newtheorem{lemma}[theorem]{Lemma}
\newtheorem{proposition}[theorem]{Proposition}
\newtheorem{corollary}[theorem]{Corollary}

\theoremstyle{definition}
\newtheorem{definition}[theorem]{Definition}
\newtheorem{example}[theorem]{Example}

\newtheorem{claim}[theorem]{Claim}
\newtheorem{question}[theorem]{Question}

\theoremstyle{remark}
\newtheorem{remark}[theorem]{Remark}

\numberwithin{equation}{section}

\newcommand{\reqnomode}{\tagsleft@false\let\veqno\@@eqno}

\begin{document}

\title{Morse resolutions and applications to lcm powers}

\author{Margherita Barile}
\address{Dipartimento di Matematica, Università degli Studi di Bari Aldo Moro, Bari, Italy
}
\email{margherita.barile@uniba.it}

\author{Trung Chau}
\address{Department of Mathematics, University of Manitoba, 420 Machray Hall, 186 Dysart
Road, Winnipeg, MB R3T 2N2 Canada}
\email{chauchitrung1996@gmail.com}

\keywords{}

\subjclass[2020]{}

\begin{abstract}
    We consider the $r$-th lcm power $I^{\{r\}}$ of a monomial ideal $I$, defined to be generated by the least common multiples of all $r$-element subsets of the set of minimal monomial generators of $I$. Some homological invariants of these ideals are studied. Specifically, we use discrete Morse theory to give  sharp upper bounds for the projective dimension and total Betti numbers of these ideals. Moreover, we completely determine the extent to which lcm powers define the homological shift ideals of edge ideals. Finally, we study when the minimal resolution of some ordinary powers and lcm powers are supported on the same CW-complex. 
\end{abstract}
\maketitle

\section{Introduction}

The study of free resolutions of monomial ideals is an important area of homological commutative algebra, in which algebraic invariants can often be investigated through combinatorial methods. In particular, the projective dimension and the Betti numbers of a monomial ideal reflect the structure of its syzygies and provide fundamental measures of its homological complexity. Although the classical Taylor resolution \cite{T} gives a universal construction, its lack of minimality motivates the search for smaller resolutions that retain a combinatorial description. In this context, discrete Morse theory has proved to be a particularly effective tool, allowing the construction of cellular resolutions by means of suitable acyclic matchings, as developed by Batzies and Welker \cite{BW02}.

In this paper, we investigate the homological properties of \textit{lcm powers} of monomial ideals. Given a monomial ideal \(I\) minimally generated by \(m_1,\ldots,m_q\), its \(r\)-th lcm power \(I^{\{r\}}\) is generated by the least common multiples of all \(r\)-element subsets of these generators. This construction is naturally related to the LCM lattice of \(I\), whose significance in the study of multigraded free resolutions was established by Gasharov, Peeva and Welker \cite{GPW99}. Our main objective is to obtain sharp upper bounds for the projective dimension and the Betti numbers of lcm powers of any monomial ideal.

\begin{theorem}[Theorem~\ref{thm:extremal-Betti}]\label{thm:extremal-Betti}
    Let $I$ be an ideal generated by $q$ monomials in a polynomial ring $R$. Then
    \[
    \pdim(R/I^{\{r\}})\leq q+1-r
    \]
    and
    \[
    \beta_i(R/I^{\{r\}})\leq \binom{r+i-2}{i-1} \binom{q}{r+i-1}
    \]
    for any integers $i,r\geq 1$. Moreover, all equalities occur for any $q,r, i\geq 1$ when $I=(x_1,\dots, x_q)$.
\end{theorem}

Our approach relies on a transfer principle for Morse resolutions, which allows suitable acyclic matchings to be transported between monomial ideals through maps preserving least common multiples. By applying this principle to squarefree Veronese ideals, we obtain sharp bounds that depend only on the number of minimal generators and the exponent of the lcm power. This provides a uniform way of controlling their homological invariants, independently of the particular structure of the original ideal. This framework, in a way, is similar to that of extremal ideals \cite{3rd,chau2026symbolicpowersintegralclosures,Lr,extremal}.

The investigation also leads to a comparison between lcm powers $I^{r+1}$ and homological shift ideals $\HS_r(I,\Bbbk)$ (see (\ref{def:shiftideals}) for the definition), which encode the multigraded shifts occurring in minimal free resolutions and have recently attracted attention in connection with the syzygies of monomial ideals \cite{Bayati2018,Bayati2023,BayatiJahaniTaghipour2019,ChauDasMaithani2026,Ficarra2025,FicarraHerzog2023,HS1,HMRZ2021}. This is motivated by the fact that we always have $\HS_r(I,\Bbbk)=I^{\{r+1\}}$ for $r\in \{0,1\}$. We remark that $\HS_r(I,\Bbbk)$ in general is dependent on the base field $\Bbbk$ of the polynomial ring $I$ is defined in, even when $I$ is an edge ideal \cite{Katzman2006}. Our next main result is that, rather surprisingly, the extent to which lcm powers determine these shifts in the case of edge ideals of graphs can be completely determined combinatorially, and in particular, is independent of the base field $\Bbbk$. 

\begin{theorem}[Theorem \ref{thm:omega-graphs}]
    Let $G$ be a connected graph with at least one edge. Then one of the following holds:
    \begin{enumerate}
        \item $G$ has an induced cycle $C_n$ for some $n\geq 3$ where $n\neq 4$ or an induced path $P_4$.
        \item $G=K_{1.n}$ for some $n\geq 1$, i.e., $G$ is the star graph on $n+1$ vertices.
        \item $G=K_{m,n}$ for some $m,n\geq 2$, i.e., $G$ is  the complete bipartite graph on two disjoint sets of $m$ and $n$ vertices.  
    \end{enumerate}
    Furthermore, the maximum index $i$ (at most the projective dimension of the edge ideal $I(G)$) such that $\HS_k(I(G),\Bbbk)=I(G)^{\{k+1\}}$ for all $k\leq i$  is $1$ if (1) holds, $n-1$ if (2) holds, and $2$ if (3) holds.

\end{theorem}

Another natural problem concerns the cellular structure of resolutions: certain ordinary powers and lcm powers have the same total Betti numbers, as follows from the results of Eliahou and Kervaire \cite{EliKer} and our computations, raising the question of whether their minimal resolutions can be supported on the same CW-complex. Our methods give a partial answer.

\begin{theorem}[Theorem~\ref{thm:q=2-or-r=2}]
    The minimal resolutions of the two ideals $(x_1,\dots, x_q)^r$ and $(x_1,\dots, x_{q+r-1})^{\{r\}}$ are supported on the same CW complex if either $q=2$ or $r=2$.
\end{theorem}

The paper is organized as follows. Section 2 develops the necessary results on Morse resolutions. Section 3 establishes the main bounds for lcm powers. Section 4 investigates their relationship with homological shift ideals, including applications to edge ideals of graphs. Section 5 addresses the question of common cellular supports for minimal resolutions.

\section*{Acknowledgements}

Chau gratefully acknowledges that this research was supported in part by the Pacific Institute for the Mathematical Sciences and a grant from the Simons Foundation International [SFI-MPS-T-Institutes-00020822, OY]. 

\textbf{AI Disclosure:} Theorem~\ref{thm:impossible-lcm-map} was proved by GPT 5.6. We take full responsibility for the content of this article.
 
\section{Preliminaries}

Throughout the section, let $S=\Bbbk[x_1,\dots, x_n]$ be  a polynomial ring over a field $\Bbbk$. Given a module $M$, a \emph{free resolution} of $M$ is a chain complex of free $S$-modules
\[
\mathcal{F}: \cdots\to F_k\xrightarrow{\partial_k}\cdots \to F_1\xrightarrow{\partial_1} F_0\xrightarrow{\partial_0} 0
\]
such that $H_i(\mathcal{F})=0$ for any $i>0$ and $H_0(\mathcal{F})\cong M$. If $\partial(F_i)\subseteq (x_1,\dots, x_n)F_{i-1}$ for any $i\geq 1$, we say that $\mathcal{F}$ is minimal, in which case the number $\beta_i^S(M)\coloneqq \operatorname{rank} F_i$ is called the \emph{$i$-th Betti number} of $M$. The subscript $S$ shall not be written if the base ring is clear. The \emph{projective dimension} of $M$, denoted by $\pdim M$, is the length of the minimal resolution of $M$. If $M$ is $\mathbb{Z}^n$-graded, then its minimal free resolution can be made $\mathbb{Z}^n$-graded, too.

Let $L=\{m_1,\dots, m_q\}$ be an ordered set of monomials. Let $\sigma=\{m_{i_1},\dots, m_{i_t}\}$ where $i_1<\cdots <i_t$ and $\sigma'\subset\sigma$ be such that $|\sigma'|=|\sigma|-1$. Then $\sigma'=\sigma\setminus \{m_{i_j}\}$ for some $j\in [t]$. We then define
\[
[\sigma\colon \sigma']\coloneqq (-1)^{j+1}.
\]
For each $i\in [q]$, let $T_i$ denote the free $S$-module with a basis $\{e_\sigma\colon \sigma\in 2^L, |\sigma|=i\}$. Let $T$ denote the sequence of maps
\[
0\to T_q\to T_{q-1}\to \cdots \to T_1\to T_0\to 0
\]
where the differentials are defined to be $\partial_0=0$ and 
\[
\partial(e_\sigma)= \sum_{\substack{\sigma'\subseteq \sigma,\\ |\sigma'|=|\sigma|-1}} [\sigma\colon \sigma'] \frac{\lcm(\sigma)}{\lcm(\sigma')} e_{\sigma'} 
\]
for any nonempty $\sigma\in 2^{L}$. Taylor \cite{T} proved that $T$ is a free resolution of $S/(L)$, thus is commonly known as a \emph{Taylor resolution} of $L$, denoted by $\TT_L$. We note that traditionally, the Taylor resolution is defined based on the minimal generating set of a monomial ideal. For our purposes, we do not always know whether a given set of monomials is a minimal generating set or not, thus we will use our unconventional version.

Given a set of monomials $L$. The power set $2^L$ can be associated with the full $(|L|-1)$-simplex. We remark that the faces of this $(|L|-1)$-simplex are in one-to-one correspondence with the power set $2^L$. Thus, any subset $\Gamma \subseteq 2^{L}$ that is closed under taking subsets corresponds to a subcomplex of the $(|L|-1)$-simplex. Such a subset $\Gamma$ is called a \emph{simplicial complex} (labeled with $L$). Note that our definition of $\Gamma$ is dependent on $L$, and a simplicial complex labeled with $L$ will be denoted by $\Gamma(L)$. When $L$ is clear from context, we will simply use $\Gamma$.

Any simplicial complex $\Gamma$ labeled with a set of monomials $L$ corresponds to a subcomplex of the Taylor resolution $\mathbb{T}_L$, which we denote by $\mathcal{F}_\Gamma$. In particular, $\TT_L=\mathcal{F}_{2^L}$, i.e., $\TT_L$ is induced by the full $(|L|-1)$-simplex. We say that $\Gamma$ \emph{supports a free resolution} of the ideal $(L)$ if $\mathcal{F}_\Gamma$ is a free resolution of $(L)$.

We recall discrete Morse theory. Again let  $L=\{m_1,\dots, m_q\}$ be a set of monomials and $\Gamma \subseteq 2^{L}$ be a simplicial complex. Let $G_{\Gamma}$ denote the directed graph with the following vertex set and edge set:
\begin{align*}
    V(G_{\Gamma})&= \Gamma,\\
    E(G_{\Gamma})&= \{\sigma\to \sigma'\colon \sigma'\subset \sigma \text{ and } |\sigma'|=|\sigma|-1\}.
\end{align*}
A set of directed edges $A\subseteq E(G_\Gamma)$ is called a \emph{homogeneous acyclic matching} of $G_\Gamma$ if
\begin{enumerate}
    \item (\emph{matching}) no two directed edges in $A$ share a vertex;
    \item (\emph{homogenenuity}) for any $(\sigma\to \sigma')\in A$, we have $\lcm(\sigma)=\lcm(\sigma')$;
    \item (\emph{acyclicity}) $G_\Gamma^A$, the directed graph $G_\Gamma$ with directed edges in $A$ being reversed, is acyclic, i.e., it does not have any directed cycle.
\end{enumerate}

It is an important theorem in discrete Morse theory that homogeneous acyclic matchings induce free resolutions. Recall that for a homogeneous acyclic matching $A$, a face $\sigma\in \Gamma$ is called \emph{$A$-critical} if it is neither the source nor the target of any edge in $A$.

\begin{theorem}[\protect{\cite[Theorem~1.3, Lemma 7.7]{BW02}}]\label{thm:Morse-theory}
    Let $L$ be a set of monomials and $\Gamma=\Gamma(L)$ a simplicial complex labeled with $L$. Assume that $\mathcal{F}_\Gamma$ is a free resolution of $S/(L)$. Then for any homogeneous acyclic matching $A$ of $G_\Gamma$, there exists a free resolution $\mathcal{F}_{\Gamma,A}$ of $S/L$ such that $(\mathcal{F}_{\Gamma,A})_i$ is a free $S$-module of rank equal to the number of $A$-critical faces of $\Gamma$. 
\end{theorem}

The differentials of $\mathcal{F}_{\Gamma,A}$ can be explicitly described (see \cite[Lemma 7.7]{BW02}). For our purposes, we shall not need this description.

Our first goal is to show that in certain circumstances, homogeneous acyclic matchings of $G_\Gamma$ can be naturally transferred to one of $G_{\Gamma'}$, where $\Gamma'$ is a simplicial complex w.r.t some other list of monomials $L'$. 

\begin{definition}
    Assume that $L=\{m_1,\dots, m_q\}$ and $L'=\{m_1',\dots, m_q'\}$ are two sets of monomials. Let $\psi:2^L\to 2^{L'}$ be the map where $\psi(\{m_{i_1},\dots, m_{i_t}\})=\{m_{i_1}
    ',\dots, m_{i_t}'\}$ for any $1\leq i_1<\cdots <i_t\leq q$.   
    We say that $\psi$ is \emph{lcm-preserving} if for any $\sigma,\tau\subseteq L$ such that $\lcm(\sigma)=\lcm(\tau)$, we have $\lcm({\psi}(\sigma))=\lcm({\psi}(\tau))$.
\end{definition}

\begin{lemma}\label{lem:Morse-transfer}
    Let $L=\{m_1,\dots, m_q\}$ and $L'=\{m_1',\dots, m_q'\}$ be two  sets of monomials, and $\Gamma(L)$ a simplicial complex labeled with $L$. Assume that the map $\psi:2^L\to 2^{L'}$, where $\psi(\{m_{i_1},\dots, m_{i_t}\})=\{m_{i_1}
    ',\dots, m_{i_t}'\}$ $1\leq i_1<\cdots <i_t\leq q$, is lcm-preserving.
    If $\Omega$ is a Morse matching for $G_{\Gamma(L)}$, then \[
    \psi(\Omega)\coloneqq \{ \psi(\sigma)\to \psi(\tau) \colon \sigma\to \tau\in \Omega\}\]
    is a Morse matching for $G_{\psi(\Gamma(L))}$. Consequently, if $\mathcal{F}_{\Gamma(L),\Omega}$ is a free resolution of $(L)$, then $\mathcal{F}_{\Gamma(L'),\psi(\Omega)}$ is a free resolution of $(\psi(L))=(L')$.
\end{lemma}

\begin{proof}
    First we show that $\psi(\Omega)$ is homogeneous. Indeed, consider $(\sigma\to  \tau)\in \Omega$. Then $\lcm(\sigma)=\lcm(\tau)$ as $\Omega$ is homogeneous. Due to our hypotheses on $\psi$, we have $\lcm(\psi(\sigma)) =  \lcm(\psi(\tau))$. In other words, $\psi(\Omega)$ is homogeneous, as claimed. Next, $\psi(\Omega)$ being a matching follows immediately from the assumption that $\Omega$ is.

    We now use contraposition to prove acyclicity. Suppose that $\psi(\Omega)$  is not acyclic. Then by \cite[Lemma~3.2.4]{Batzies2002}, there exists a directed cycle
    \[
    \sigma_1' \to \tau_1' \ {\color{blue} \to}\  \sigma_2' \to \tau_2' \ {\color{blue} \to}\  \cdots \ {\color{blue} \to}\  \sigma_l'\to \tau_l' \ {\color{blue} \to}\ \sigma_1'
    \]
    for some $l\geq 1$, where $(\sigma_{i+1}'\to \tau_i') \in \psi(\Omega)$, and $(\sigma_i'\to \tau_i') \in E(G_{\psi(\Gamma(L))})\setminus \psi(\Omega)$,  for any $i\in [l]$. As $\psi$ is a bijection, there exist preimages $\sigma_i$ and $\tau_i$ such that $\psi(\sigma_i)=\sigma_i'$ and $\psi(\tau_i)=\tau_i'$ for any $i\in [l]$. This implies that $ (\sigma_{i+1} \ {\color{blue} \to}\  \tau_{i} ) \in \Omega$ for any $i\in [l]$, where we adopt the convention $l+1=1$. Since $\psi$ preserves inclusion, we obtain $\sigma_i \to \tau_i \in E(G_{\Gamma(L)})\setminus \Omega$. Putting all these directed edges together we obtain a directed cycle in $G_{\Gamma(L)}^{\Omega}$. In other words, $\Omega$ is not acyclic, as~desired.
\end{proof}

\begin{example}
    Consider the algebra map that extends from the following:
    \begin{align*}
        \psi\colon R=\Bbbk[w,x,y,z]&\to S=\Bbbk[a,b,c,d]\\
        w&\mapsto ab\\
        x&\mapsto c\\
        y&\mapsto d\\
        z&\mapsto 1.
    \end{align*}
    Now consider the set $L=\{wx,wz,xy,yz\}$. The map $\psi$ induces a natural map $2^L\to 2^{\psi(L)}$, which we will also denote by $\psi$. The hypothesis of \cref{lem:Morse-transfer} can be easily verified. It is known from \cite[Example 2.13]{CK24} that
    \begin{multline*}
        \Omega = \{ \{wx,wz,xy,yz\}\to \{wx,wz,xy\}, \quad \{wx,xy,yz\}\to \{wx,yz\}, \quad \{wz,xy,yz\}\to \{wz,xy\} \}
    \end{multline*}
    is a homogeneous acyclic matching of $G_L$. By our preceding result, the set
    \[
    \psi(\Omega) = \{ \{abc,ab,cd,d\}\to \{abc,ab,cd\}, \quad \{abc,cd,d\}\to \{abc,d\}, \quad \{ab,cd,d\}\to \{ab,cd\} \}
    \]
    is a homogeneous acyclic matching of $G_{\{abc,ab,cd,d\}}$. The Morse resolutions induced from $\Omega$ and $\psi(\Omega)$ are
    \[
    R\to R^4\to R^4\to R\to 0
    \]
    and 
    \[
    S\to S^4\to S^4\to S\to 0,
    \]
    respectively. However, only the first one is minimal (\cite[Example 2.13]{CK24}). The second is not minimal simply from the observation that the ideal generated by $\psi(L)=\{abc,ab,cd,d\}$ is \emph{minimally} generated by only two monomials, namely $ab$ and $d$.
\end{example}

\begin{corollary}\label{cor:betti-bound-Morse-transfer}
    Assume the same set-up as in Lemma~\ref{lem:Morse-transfer}. Further assume that $(L)$ has a minimal Morse resolution. Then
    \[
    \beta_i((\psi(L)))\leq \beta_{i}((L))
    \]
    for any $i\in \mathbb{Z}$.
\end{corollary}

\begin{remark}\label{rem:lcm-lattice}
    For a monomial ideal $I$, let $\mingens(I)$ denote the set of minimal monomial generators of $I$, and $\operatorname{LCM}(I)$ denote the lattice with elements $\{ \lcm(\sigma)\colon \sigma\in 2^{\mingens(I)} \}$ and the join and meet operations being taking the lcm and gcd, respectively. This is called the \emph{LCM lattice} of $I$. Gasharov, Peeva, and Welker \cite{GPW99} showed that the LCM lattice of $I$ determines all of its Betti numbers. From the perspective of LCM lattices, if $L$ is a list of monomials and $\psi$ is a map as in the hypothesis of \cref{lem:Morse-transfer}, \cite[Theorem 3.3]{GPW99} implies Corollary~\ref{cor:betti-bound-Morse-transfer}. Lemma~\ref{lem:Morse-transfer} complements this fact, stating that a Morse resolution of $(L)$ is a Morse resolution of $(\psi(L))$, after a change of rings. 
\end{remark}

Another situation where Morse resolutions are ``transfered" is in the case of restriction ideals. Recall that for a monomial ideal $I$ and a monomial $m$, the \emph{restriction} of $I$ with respect to $m$, denoted by $I^{\leq m}$, is the ideal generated by all minimal monomial generators of $I$ that divide $m$.

The following is likely known to experts. We provide a proof for completion.

\begin{lemma}\label{lem:Morse-restriction}
    If a monomial ideal $I$ has a minimal Morse resolution, then so does $I^{\leq m}$ for any monomial $m$.
\end{lemma}

\begin{proof}
    Assume that $I$ has a minimal Morse resolution. By definition, there exists a simplicial complex $\Gamma(I)$ labeled with minimal monomial generators of $I$, and a Morse matching $\Omega$ such that $\mathcal{F}_{\Gamma(L),\Omega}$ is a minimal Morse resolution of $I$. We shall use these to construct a minimal Morse resolution of $I^{\leq m}$.
    
    Let $\Gamma(I)^{\leq m}$ denote the subcomplex of $\Gamma(I)$ containing only the faces whose labels divide $m$. By definition, $\Gamma(I)^{\leq m}$ is a simplicial complex labeled with minimal monomial generators of $I^{\leq m}$. Moreover, by \cite[Lemma 4.4]{HHZ-Dirac}, $\Gamma(I)^{\leq m}$ supports a free resolution of $I^{\leq m}$. Next we define
    \[
    \Omega^{\leq m} = \{\sigma \to \sigma'\in \Omega\colon \lcm(\sigma)\mid m \}.
    \]
    Since $\Omega$ is a homogeneous matching, so is any of its subset, and thus $\Omega^{\leq m}$ in particular. Moreover, it is routine to verify that $G_{\Gamma(L)^{\leq m}}^{\Omega^{\leq m}}$ is a directed subgraph of $G_{\Gamma(L)}^{\Omega}$. Hence the former is acyclic as the latter is. Thus $\Omega^{\leq m}$ is a homogeneous acyclic matching of $I^{\leq m}$. It now suffices to show that $\mathcal{F}_{\Gamma(I)^{\leq m},\Omega^{\leq m}}$ is a minimal resolution. Indeed, the ranks of the free modules in $\mathcal{F}_{\Gamma(I)^{\leq m},\Omega^{\leq m}}$, by our construction, are exactly the numbers of critical faces of $\Gamma(I)$ whose labels divide $m$, which are exactly the Betti numbers of $I^{\leq m}$ by \cite[Lemma 4.4]{HHZ-Dirac}. This concludes the~proof.
\end{proof}

\section{Applications to lcm powers of monomial ideals}

From a certain perspective, \cref{lem:Morse-transfer} allows us to construct a Morse resolution of $R/J$ for a monomial ideal $J$ of a polynomial ring $R$, given a monomial ideal $I$ of a polynomial ring $S$, a Morse resolution for $S/I$, and a map $\psi\colon 2^{\mingens(I)}\to 2^{M}$ where $M$ is a monomial generating set of $J$. 

We present an application to a newly constructed class of ideals. For the rest of the section let $I$ denote an ideal minimally generated by $m_1,\dots, m_q$ for some integer $q\geq 1$, and $1\leq r\leq q$ denote an integer. 

We  define the \textit{lcm $r$-th power} of $I$ to be
\[
I^{\{r\}}\coloneqq (\lcm(m_{i_1},\dots, m_{i_r}) \colon 1\leq i_1<\cdots <i_r\leq q). 
\]
We remark that $I^{\{1\}}=I$. Moreover, set $\psi_{I,r}$ to be the unique map
\begin{align*}
    \psi_{I,r}:2^{ \mingens((x_1,\dots, x_q)^{\{r\}}) } \to 2^{ \{\lcm(m_{i_1},\dots, m_{i_r})\colon 1\leq i_1<i_2<\cdots < i_r\leq q\} }
\end{align*}
where for each $\sigma\in 2^{ \mingens((x_1,\dots, x_q)^{\{r\}})}$, we have \[
x_{i_1}\cdots x_{i_r} \in \sigma \quad \text{if and only if}  \quad\lcm(m_{i_1},\dots, m_{i_r})\in \psi_{I,r}(\sigma).
\]
\begin{lemma}\label{lem:psi-I-r-lcm-preserve}
    The map $\psi_{I,r}$ preserves least common multiples.
\end{lemma}

\begin{proof}
    Consider two subsets $\sigma,\tau$ of minimal generators of $(x_1,\dots, x_q)^{\{r\}}$ such that $\lcm(\sigma)=\lcm(\tau)$. By definition,  there exist two collections $\mathcal{I}_\sigma$ and $\mathcal{I}_\tau$ of subsets of $[q]$ of size $r$ such that 
    \[
    \sigma = \{ \prod_{i\in \mathcal{I}} x_i \mid \mathcal{I} \in \mathcal{I}_\sigma \} \quad \text{and} \quad \tau = \{ \prod_{i\in \mathcal{I}} x_i \mid \mathcal{I} \in \mathcal{I}_\tau \}. 
    \]
    Then 
    \[
    \lcm(\sigma) = \prod_{i\in \cup_{\mathcal{I}\in \mathcal{I}_\sigma} \mathcal{I{}}} x_i \quad \text{and} \quad \lcm(\tau) = \prod_{i\in \cup_{\mathcal{I}\in \mathcal{I}_\tau} \mathcal{I{}}} x_i. 
    \]
    Thus $\lcm(\sigma)=\lcm(\tau)$ simply means that $\bigcup_{\mathcal{I}\in \mathcal{I}_\sigma} \mathcal{I} = \bigcup_{\mathcal{I}\in \mathcal{I}_\tau} \mathcal{I}$. Therefore,
    \begin{align*}
        \lcm(\psi_{I,r}(\sigma)) &= \lcm(\{ \lcm(\{m_i\mid i\in \mathcal{I}\}) \mid \mathcal{I}\in \mathcal{I}_\sigma  \})\\
        &= \lcm(\{m_i\mid i\in \cup_{\mathcal{I\in \mathcal{I}_\sigma}} \mathcal{I}\})\\
        &= \lcm(\{m_i\mid i\in \cup_{\mathcal{I\in \mathcal{I}_\tau}} \mathcal{I}\})\\
        &= \lcm(\{ \lcm(\{m_i\mid i\in \mathcal{I}\}) \mid \mathcal{I}\in \mathcal{I}_\tau  \})\\
        &=\lcm(\psi_{I,r}(\tau)),
    \end{align*}
    as desired.
\end{proof}

We have obtained the condition to apply Lemma~\ref{lem:Morse-transfer} to $\psi_{I,r}$. In fact, we know more about the ``source" ideal $(x_1,\dots, x_q)^{\{r\}}$. Recall that a homogeneous ideal $I$ \emph{has a linear resolution} if in the minimal free resolution of $I$, all the entries in the representation matrices are linear. We also recall the ``restriction lemma".

\begin{lemma}[\protect{\cite[Lemma~4.4]{HHZ-Dirac}}]\label{lem:restriction}
    Let $I$ be a monomial ideal and $m$ a monomial  in $\Bbbk[x_1,\dots, x_n]$. Let $\mathcal{F}$ be a $\mathbb{Z}^n$-graded minimal free resolution of $I$, and $\mathcal{F}^{\leq m}$ be the subcomplex of $\mathcal{F}$ consisting of all basis elements with $\mathbb{Z}^n$-degree at most the exponent vector of $m$. Then $\mathcal{F}^{\leq m}$ is a minimal free resolution of $I^{\leq m}$.
\end{lemma}

\begin{lemma}\label{lem:Morse-lcm-power-maximal-ideal}
    For any two integers $q,r\geq 1$, the ideal $(x_1,\dots, x_q)^{\{r\}}$ has a linear minimal Morse resolution.
\end{lemma}

\begin{proof}
    It is straightforward that $(x_1,\dots, x_q)^{\{r\}}=\left((x_1,\dots, x_q)^{r}\right)^{\leq \prod_{i=1}^q x_i}$. We know that $(x_1,\dots, x_q)^{r}$ has a linear minimal Morse resolution by \cite[Fact 5.2]{BW02} and \cite{EliKer}. Thus so does $(x_1,\dots, x_q)^{\{r\}}$ by Lemmas~\ref{lem:Morse-restriction} and \ref{lem:restriction}.
\end{proof}

Recall that a decomposition $I=J+K$ of a monomial ideal $I$ is called a \emph{Betti splitting} if 
\[
\beta_{i}(I) = \beta_{i}(J) + \beta_{i}(K) + \beta_{i-1}(J\cap K)
\]
for any integer $i$. We recall a sufficient condition to obtain Betti splittings.

\begin{lemma}[\protect{\cite[Corollary~2.7]{francisco-2009}}]\label{lem:betti-splitting}
    Let $I=xJ+K$ where $I$ is a monomial ideal, $x$ is a variable, and $J$ and $K$ are monomial ideals whose minimal monomial generators do not involve $x$, such that $\mingens(I)=\mingens(xJ)\sqcup \mingens(K)$. If $J$ has a linear resolution, then $I=xJ+K$ is a Betti splitting.
\end{lemma}

We are in a position to prove the main result of this section.

\begin{theorem}\label{thm:extremal-Betti}
    Let $I$ be a monomial ideal generated by $q$ monomials in a polynomial ring $R$. Then
    \[
    \pdim(R/I^{\{r\}})\leq q+1-r
    \]
    and
    \[
    \beta_i(R/I^{\{r\}})\leq \binom{r+i-2}{i-1} \binom{q}{r+i-1}
    \]
    for any integers $i\geq 1$ and $r\geq 1$. Moreover, equalities occur for each $q,r,$ and $i\geq 1$ when $I=(x_1,\dots, x_q)$.
\end{theorem}

\begin{proof}
    By Lemmas~\ref{lem:psi-I-r-lcm-preserve} and \ref{lem:Morse-lcm-power-maximal-ideal} and Corollary~\ref{cor:betti-bound-Morse-transfer}, it suffices to prove the equality in the case where $I=(x_1,\dots, x_q)$ in the ring $R=\Bbbk[x_1,\dots, x_q]$. Set $\beta(i,r,q)\coloneqq \beta_i(R/I^{\{r\}})$. We will find a recursive formula for this sequence.
    
    Set
    \begin{align*}
        J&\coloneqq (x_1,\dots, x_{q-1})^{\{r-1\}},\\
        K&\coloneqq (x_1,\dots, x_{q-1})^{\{r\}}.
    \end{align*}
    It is clear that $\mingens(I^{\{r\}})=\mingens(x_qJ)\sqcup \mingens(K)$ and $x_qJ\cap K= x_q (x_1,\dots, x_{q-1})^{\{r\}}$.   
    As $J$ has a linear resolution by Lemma~\ref{lem:Morse-lcm-power-maximal-ideal}, the decomposition $I^{\{r\}}=x_qJ+K$ is a Betti splitting by Lemma~\ref{lem:betti-splitting}. For any $r,q,i\geq 1$, we thus have
    \begin{align*}
        \beta(i,r,q) = \beta_i(R/I^{\{r\}}) &= \beta_i(R/x_qJ) + \beta_i(R/K) + \beta_{i-1}(R/(x_qJ\cap K))\\
        &=  \beta(i,r-1,q-1) + \beta(i,r,q-1)  + \beta(i-1,r,q-1).
    \end{align*}
    Now that we have a recursive formula, we can prove our result with induction on $i,q,$ and $r$. First we consider the base case. When $i=1$, it is well-known that $\beta(i,r,q)$ is simply the minimal number of generators of $I^{\{r\}}$, i.e.,
    \[
    \beta(1,r,q) = \binom{q}{r} = \binom{r+1-2}{r-1}\binom{q}{r+1-1},
    \]
    as desired. When $r=1$, the ideal $I^{\{r\}}=(x_1,\dots,x_q)$ has a minimal Taylor resolution, and thus
    \[
    \beta(i,1,q) = \binom{q}{i} = \binom{1+i-2}{i-1} \binom{q}{1+i-1},
    \]
    as desired. When $q=1$, for the ideal $I^{\{r\}}$ to be nonzero, we need $r=1$ as well. The result then follows from the previous case. Now we can compute $\beta(i,r,q)$ using the induction hypothesis:
    \begin{align*}
        \beta(i,r,q)&=\beta(i,r-1,q-1) + \beta(i,r,q-1)  + \beta(i-1,r,q-1) \\
        &=  \binom{r+i-3}{i-1}\binom{q-1}{r+i-2} + \binom{r+i-2}{i-1}\binom{q-1}{r+i-1} + \binom{r+i-3}{i-2}\binom{q-1}{r+i-2}.
    \end{align*}
    Adding the first and third terms, which we observe to have an obvious common factor, and using $\binom{n+1}{k+1}=\binom{n}{k+1}+\binom{n}{k}$, we obtain
    \begin{align*}
         \beta(i,r,q) =  \binom{r+i-2}{i-1}\binom{q-1}{r+i-2} + \binom{r+i-2}{i-1}\binom{q-1}{r+i-1}.
    \end{align*}
    Using the formula $\binom{n+1}{k+1}=\binom{n}{k+1}+\binom{n}{k}$ again, we obtain
    \begin{align*}
         \beta(i,r,q) =   \binom{r+i-2}{i-1}\binom{q}{r+i-1},
    \end{align*}
    as desired. The formula for projective dimension follows immediately from those for Betti numbers. This concludes the proof.
\end{proof}

As an application, we determine a sharp upper bound for the first homological shift ideal. We first recall  the concept. Let $I$ be a monomial ideal in $R=\Bbbk[x_1,\dots, x_n]$. Then it is known that the minimal free resolution
\[
\mathcal{F}:\cdots \to F_1\to F_0
\]
is $\mathbb{Z}^n$-graded, i.e., we can  uniquely write
\[
F_i=\bigoplus_{\mathbf{a}\in \mathbb{Z}^n} R(-\mathbf{a})^{\beta_{i,\mathbf{a}}(I)},
\]
for any $i\geq 0$. For each $i\geq 0$, we can then define
\begin{equation}\label{def:shiftideals}
    \HS_i(I,\Bbbk)\coloneqq ( \mathbf{x}^{\mathbf{a}} \colon \beta_{i,\mathbf{a}}(I)>0 ).
\end{equation}
It is immediate that $\HS_0(I,\Bbbk)=I$ itself. It is widely accepted that homological shift ideals are hard to compute, as the only known method is to construct the minimal free resolution first and extract its information. Nevertheless, the formula for $\HS_1(I,\Bbbk)$ is known. We then obtain sharp bounds on its Betti numbers and projective dimension.

\begin{lemma}[\protect{\cite[Lemma~2.2]{HS1}}]\label{lem:HS1}
    For any monomial ideal $I$, we have $\HS_1(I,\Bbbk)=I^{\{2\}}$.
\end{lemma}

\begin{corollary}
    For any monomial ideal $I$ with $q$ generators, we have
    \[
    \pdim(\HS_1(I,\Bbbk))\leq q-2
    \]
    and
    \[
    \beta_i^R(\HS_1(I,\Bbbk))\leq (i+1) \binom{q}{i+2}
    \]
    for any integers $i\geq 1$ and $r\geq 1$. Moreover, equalities occur for each $q,r,$ and $i\geq 1$ when $I=(x_1,\dots, x_q)$.
\end{corollary}
\begin{proof}
    The result follows from Theorem~\ref{thm:extremal-Betti} and Lemma~\ref{lem:HS1}. 
\end{proof}

\section{Lcm powers versus homological shift ideals}

Since the minimal free resolution of a monomial ideal $I$ (of a polynomial ring over a field $\Bbbk$) is a direct summand of its Taylor resolution, we always have $\HS_i(I,\Bbbk)\subseteq I^{\{i+1\}}$ for any integer $i$. In this section, we shall study when this equality occurs. Throughout  this section, unless stated otherwise, $I$ denotes a monomial ideal of a polynomial ring over a field $\Bbbk$. We recall a  key lemma, which follows directly from Lemma~\ref{lem:restriction}.

\begin{lemma}\label{lem:restriction-lemma}
    For any integer $i$, monomial ideal $I$, and monomial $m$, we have
    $\HS_i(I,\Bbbk)^{\leq m}=\HS_i(I^{\leq m},\Bbbk)$
\end{lemma}

It is straightforward to characterize when we always have an equality. For a monomial ideal $I$, let $\mu(I)$ denote the number of minimal generators of $I$. Note that we already have $\HS_i(I,\Bbbk)=I^{\{i+1\}}$ for $i=0,1$ and we always have $\HS_i(I,\Bbbk)=I^{\{i+1\}}=0$ for any $i\geq \mu(I)$.

\begin{proposition}\label{prop:Taylor-equality-lcm-power}
    For a monomial ideal $I$ and an integer $i$, we have $\HS_i(I,\Bbbk)\subseteq I^{\{i+1\}}$. Moreover, the following are equivalent:
    \begin{enumerate}
        \item $\HS_i(I,\Bbbk)= I^{\{i+1\}}$ for any $i\in \mathbb{N}$.
        \item $\HS_{\mu(I)-1}(I,\Bbbk)= I^{\{\mu(I)\}}$.
        \item The Taylor resolution of $I$ is minimal.
    \end{enumerate}
\end{proposition}

\begin{proof}
    The containment follows from the description of the Taylor resolution of $I$, and the fact that the minimal resolution of $I$ is a direct summand of it.
    
    $(1)\implies (2)$: Straightforward.
    
    $(2)\implies (3)$: We shall show this with contraposition.  Assume the Taylor resolution of $I$ is not minimal. Assume that $I$ is minimally generated by monomials $m_1,\dots, m_q$. We have the following claim.
    \begin{claim}\label{clm:taylor-not-minimal}
        We have $\beta_{q-1}(I)=0$.
    \end{claim}
    \begin{proof}[Proof of Claim~\ref{clm:taylor-not-minimal}]
        Since the Taylor resolution of $I$ is not minimal, by the description of the differentials, we have $\lcm(\sigma)=\lcm(\sigma\setminus \{m_k\})$ for some $\sigma\subseteq \{m_1,\dots, m_q\}$ and some $i\in [q]$ such that $m_k\in \sigma$. By adding all the monomials not in $\sigma$ to both sets, we have $\lcm(M)=\lcm(M\setminus \{m_k\})$, where $M\coloneqq \{m_1,\dots, m_q\}$. Now set $\Omega\coloneqq \{ M\to M\setminus \{m_k\} \}$. The condition $\lcm(M)=\lcm(M\setminus \{m_k\})$ guarantees that $\Omega$ is a homogeneous, and it is straightforward that $\Omega$ is also an acyclic matching. Note that $\mathcal{F}_{2^M, \Omega}$  is of length $q-2$. Therefore, by Theorem~\ref{thm:Morse-theory}, we  have $\pdim I \leq q-2$. In particular, we have $\beta_{q-1}(I)=0$, as desired.
    \end{proof}
    By Claim~\ref{clm:taylor-not-minimal}, we have 
    \[
    \HS_{q-1}(I,\Bbbk) =0 \neq (\lcm(m_1,\dots, m_q)) = I^{\{q\}},
    \]
    as desired.

    $(3)\implies (1)$: This follows directly from definition and the description of the Taylor resolution.
\end{proof}

Next, we show that the two ideals $\HS_i(I,\Bbbk)$ and $I^{\{i+1\}}$ can coincide up until a point, and differ afterwards (up until index $\mu(I)-1$). 

\begin{lemma}\label{lem:construction-HS-lcm-power-coincide}
    Assume that there exist a monomial ideal $I$ and two integers $k$ and $q$ such that the following are satisfied:
    \begin{enumerate}
        \item $1\leq k\leq q-2$;
        \item $I$ has $q$ minimal generators;
        \item $\HS_i(I,\Bbbk) \subsetneq I^{\{i+1\}}$ if and only if $k<i<q$.
    \end{enumerate}
    Then for a variable $x$ that does not divide any minimal generator of $I$, we have $\HS_i(I+(x))\subsetneq \left(I+(x) \right)^{\{i+1\}}$ if and only if $k<i< q+1$. 
\end{lemma}

\begin{proof}
    By Lemma~\ref{lem:betti-splitting}, $I+(x)$ is a Betti splitting of the ideal $I+(x)$ with $I\cap (x)=xI$. Thus for any integer $i$, we have
    \[
    \HS_i(I+(x))= \HS_i(I,\Bbbk) + \HS_i((x),\Bbbk)+ x\HS_{i-1}(I,\Bbbk).
    \]
    Now assume that $2\leq i\leq k$, we have
    \begin{align*}
        \HS_i(I+(x))&= \HS_i(I,\Bbbk) + \HS_i((x),\Bbbk)+ x\HS_{i-1}(I,\Bbbk)\\
        &=\HS_i(I,\Bbbk)+ x\HS_{i-1}(I,\Bbbk)\\
        &=I^{\{i+1\}} +xI^{\{i\}}\\
        &=(I+(x))^{\{i+1\}}.
    \end{align*}
    On the other hand, for $k<i\leq q$ (which implies that $i\geq 1$, thus $\HS_i((x),\Bbbk)=0$), we have
    \begin{align*}
        \HS_i(I+(x))&= \HS_i(I,\Bbbk) + \HS_i((x),\Bbbk)+ x\HS_{i-1}(I,\Bbbk)\\
        &=\HS_i(I,\Bbbk)+ x\HS_{i-1}(I,\Bbbk)\\
        &\subsetneq I^{\{i+1\}} +xI^{\{i\}}\\
        &= (I+(x))^{\{i+1\}}.
    \end{align*}
    The result then follows.
\end{proof}

\begin{theorem}
    For each pair of integers $(q,k)$ where $1\leq k\leq q-1$, there exists a monomial ideal $I$ with $q$ minimal generators such that $\HS_i(I,\Bbbk)\subsetneq I^{\{i+        1\}}$ if and only if $k<i<q$.
\end{theorem}

\begin{proof}
    If $k=q-1$, the ideal $I=(x_1,\dots, x_q)$ satisfies $\HS_i(I,\Bbbk)= I^{\{i+        1\}}$ for any $i$, as desired.

    If $k=q-2$, we consider the ideal $J=(x_1\cdots x_{q-1},x_1x_q,x_2x_q,\dots, x_{q-1}x_q)$. By Lemma~\ref{lem:betti-splitting}, the decomposition $J=x_qJ_1+J_2$, where   $J_1=(x_1,\dots, x_{q-1})$ and $J_2=(x_1\cdots x_{q-1})$, is a Betti splitting with $x_qJ_1\cap J_2=x_qJ_2$. Thus for any index $i$, we have
    \[
    \HS_i(J,\Bbbk)=x_q\HS_i(J_1,\Bbbk)+\HS_{i}(J_2,\Bbbk)+x_q\HS_{i-1}(J_2,\Bbbk).
    \]
    In particular, if $i\geq 2$, then $\HS_i(J,\Bbbk)=x_q\HS_i(J_1,\Bbbk)$. For $2\leq i\leq q-2$, we then have
    \begin{align*}
        \HS_i(J,\Bbbk)= x_q\HS_{i}(J_1,\Bbbk)=x_qJ_1^{\{i+1\}} = J^{i+1},
    \end{align*}
    where the second equality is due to Proposition~\ref{prop:Taylor-equality-lcm-power} and the fact that $J_1$ has a minimal Taylor resolution (see, e.g., \cite[Corollary 4.9]{Alesandroni2017}).
    On the other hand, we have
    \[
    \HS_{q-1}(J,\Bbbk)=x_q\HS_{q-1}(J_1,\Bbbk)=0\neq (x_1\cdots x_q) = J^{\{q\} }.
    \]
    Thus $\HS_i(J,\Bbbk)\subsetneq J^{\{i+        1\}}$ if and only if $k=q-2<i<q$, as desired.

    Now we can assume that $k\leq q-3$. Let $K$ be a monomial ideal with $k+2$ minimal generators such that $\HS_i(K,\Bbbk)\subsetneq K^{\{i+        1\}}$ if and only if $k<i<k+2$. Such an ideal $K$ exists due to the previous case, and note that $k+2<q$. The ideal $K+(y_1,y_2,\dots, y_{q-k-2} )$, where none of the variables $y_1,\dots, y_{q-k-2}$ divides any minimal monomial generator of $K$, then satisfies the desired property by Lemma~\ref{lem:construction-HS-lcm-power-coincide}. This concludes the proof. 
\end{proof}

The preceding theorem in particular shows that the phenomenon $\HS_i(I,\Bbbk)=I^{ \{i+1\} }$ can occur up to any point. Conversely, given a monomial ideal $I$, what is the largest index (dependent on $I$ itself) such that up until this index, we always have the equality $\HS_i(I,\Bbbk)=I^{ \{i+1\} }$? We shall investigate this question and give a complete answer for the most popular class of monomial ideals: edge ideals of graphs.

For the rest of the section, fix a field $\Bbbk$. For a finite simple graph $G=(V(G),E(G))$, the \emph{edge ideal} of $G$, denoted by $I(G,\Bbbk)$, is defined to be the ideal
\[
\langle xy\mid \{x,y\}\in E(G) \rangle \subseteq \Bbbk[V(G)].
\]
By definition, isolated vertices do not add any edges to an edge ideal. We therefore shall always assume that our graphs do not have any isolated vertices.

\begin{definition}\label{def:omega-def}
    For a finite simple graph $G$ with at least one edge, we define
    \[
\omega(G,\Bbbk)\coloneqq  \max\{ i\leq \pdim I(G)\mid \HS_k(I(G),\Bbbk)=I(G)^{\{k+1\}}\ \forall k\leq i \}
\]
    and call it the \emph{predictability index of homological shift} of $G$.
\end{definition}

\begin{remark}
    In general, $\omega(G,\Bbbk)$ may be dependent on $\Bbbk$, which makes it not possible to compute with Macaulay2 \cite{M2}, as this program requires the specification of the base field. However, when $G$ has at most 8 vertices, the minimal resolution of $I(G)$ does not depend on $\Bbbk$ (\cite[Theorem 4.2]{ChauHaMaithani2026}), and thus can be computed with Macaulay2 \cite{M2}. In particular, $\omega(G,\Bbbk)=\omega(G,\mathbb{Q})$ can also be computed in this case. Throughout, if $G$ has at most 8 vertices, we shall only present the number $\omega(G,\Bbbk)$ with no  explanation.
\end{remark}

The goal for the rest of the section is to completely determine $\omega(G,\Bbbk)$. In fact, we have computed this invariant in a number of cases.

\begin{lemma}\label{lem:Taylor-omega}
    Let $G$ be a graph with at least one edge. Then $I(G)$ admits a minimal Taylor resolution if and only if $\omega(G,\Bbbk)=|E(G)|-1$.
\end{lemma}
\begin{proof}
    This is a special case of Proposition~\ref{prop:Taylor-equality-lcm-power}.
\end{proof}

It is standard (see, e.g., \cite[Lemma~2.1]{JacquesKatzman2005}) that if $I,J$ are monomial ideals generated in disjoint sets of variables in a polynomial ring $S$, then the minimal free resolution of $S/(I+J)$ is the tensor product of those of $S/I$ and $S/J$. Thus we obtain the following immediately.

\begin{lemma}
    Let $G$ be a graph with connected component $G_1,\dots, G_n$. Set $S=\Bbbk[V(G)]$. Then for any $i\geq 0$, we have
    \[
    \HS_{i}(S/I(G),\Bbbk) = \sum_{i_1+\cdots +i_n=i} \HS_{i_1}(S/(I(G_{1}),\Bbbk)\cdots  \HS_{i_n}(S/I(G_{n}),\Bbbk).
    \]
    In particular, we have $\omega(G,\Bbbk)=\min \{\omega(G_i,\Bbbk)\mid i=1,\dots, n  \}$.
\end{lemma}

Due to the above lemma, we shall only work on connected graphs. To compute $\omega(G,\Bbbk)$ for any connected graph $G$, we shall need new tools. The two following key lemmas give bounds for this invariant.

\begin{lemma}\label{lem:restriction-omega}
    Let $H$ be an induced subgraph of $G$. If $\pdim I(H)>\omega(H,\Bbbk)$, then $\omega(G,\Bbbk)\leq \omega(H,\Bbbk)$.
\end{lemma}

\begin{proof}
    Since $\pdim I(H)>\omega(H,\Bbbk)$, we have 
    \[\HS_{\omega(H,\Bbbk)+1}(I(G))^{\leq \prod_{x\in V(H)}x} = \HS_{\omega(H,\Bbbk)+1}(I(H))\neq I(H)^{\{\omega(H,\Bbbk)+2\}} = (I(G)^{\{\omega(H,\Bbbk)+2\}})^{\leq \prod_{x\in V(H)}x},\]
    where the two equalities are due to the restriction lemma. In particular, this implies that $\HS_{\omega(H,\Bbbk)+1}(I(G))\neq I(G)^{\{\omega(H,\Bbbk)+2\}}$. Thus $\omega(G,\Bbbk)\leq \omega(H,\Bbbk)$ by definition, as desired.
\end{proof}

\begin{lemma}\label{lem:omega-at-least-1}
    Let $G$ be a graph with at least one edge. Then $\omega(G,\Bbbk)=0$ if and only if $G$ has exactly one edge. Consequently, if $G$ has more than one edge, then $\omega(G,\Bbbk)\geq 1$.
\end{lemma}

\begin{proof}
    It is well-known that $\pdim I(G) \geq 0$ with equality occurring exactly when $G$ has exactly one edge. The result then follows from the fact that $\HS_0(I(G)) = I(G)=I(G)^{\{1\}}$ and $\HS_1(I(G))=I(G)^{\{2\}}$ for any graph $G$ (Lemma~\ref{lem:HS1}).
\end{proof}

A direct corollary of the two preceding lemmas is the following.

\begin{corollary}\label{cor:omega=1}
    Let $G$ be a graph. Assume that $G$ has more than one edge and there exists an induced subgraph $H$ of $G$ such that $\omega(H,\Bbbk)=1<\pdim I(H)$. Then $\omega(G,\Bbbk)=1$.
\end{corollary}

This corollary allows us to determine $\omega(G,\Bbbk)$ easily if $G$ contains certain induced subgraphs. We  provide a list of such graphs in our next result. Before that, we recall some standard graph names.

For each $n\geq 1$, the \emph{$n$-path} (denoted by $P_n$), the \emph{$n$-cycle} (denoted by $C_n$, only defined for $n\geq 3$), and the complete graph on $n$ vertices (denoted by $K_n$), are graphs with the vertex set $\{x_1,\dots, x_n\}$ and edge sets
\begin{align*}
    E(P_n)&=\{\{x_i,x_{i+1}\}\mid i\in [n-1]\},\\
    E(C_n)&=E(P_n)\cup \{\{x_1,x_n\}\},\\
    E(K_n)&= \{ \{x_i,x_j \}\mid 1\leq i<j\leq n \},
\end{align*}
respectively. We also recall some small graphs below.

\begin{figure}[h!]
\centering
\begin{tabular}{ccccccc}
\begin{tikzpicture}[scale=1, every node/.style={circle,draw,fill=black,inner sep=1pt}]
\node (B) at (0.9,0) {};
\node (C) at (-0.9,0) {};
\node (A) at (0,1.5) {};
\node (D) at (0,2.7) {};

\draw (A)--(B)--(C)--(A)--(D);

\end{tikzpicture}
&&
\begin{tikzpicture}[scale=1, every node/.style={circle,draw,fill=black,inner sep=1pt}]
\node (B) at (1,0) {};
\node (D) at (-1,0) {};
\node (A) at (0,1.5) {};
\node (C) at (0,-1.5) {};

\draw (D)--(A)--(B)--(C)--(D)--(B);

\end{tikzpicture}&&
\begin{tikzpicture}[scale=1, every node/.style={circle,draw,fill=black,inner sep=1pt}]
\node (A) at (0,0) {};
\node (B) at (0,-1) {};
\node (C) at (0,-2) {};
\node (D) at (0.7,0.8) {};
\node (E) at (-0.7,0.8) {};

\draw (D)--(A)--(B)--(C);
\draw (E)--(A);

\end{tikzpicture}&&
\begin{tikzpicture}[scale=1, every node/.style={circle,draw,fill=black,inner sep=1pt}]
\node (A) at (0,0) {};
\node (B) at (1.4,0) {};
\node (C) at (1.4,1.4) {};
\node (D) at (0,1.4) {};
\node (E) at (1.4,-1.4) {};

\draw (D)--(A)--(B)--(C)--(D);
\draw (E)--(B);

\end{tikzpicture}
\end{tabular}
\caption{The paw graph, the diamond graph, the fork graph, and the banner graph, from left to right, respectively.}
\end{figure}
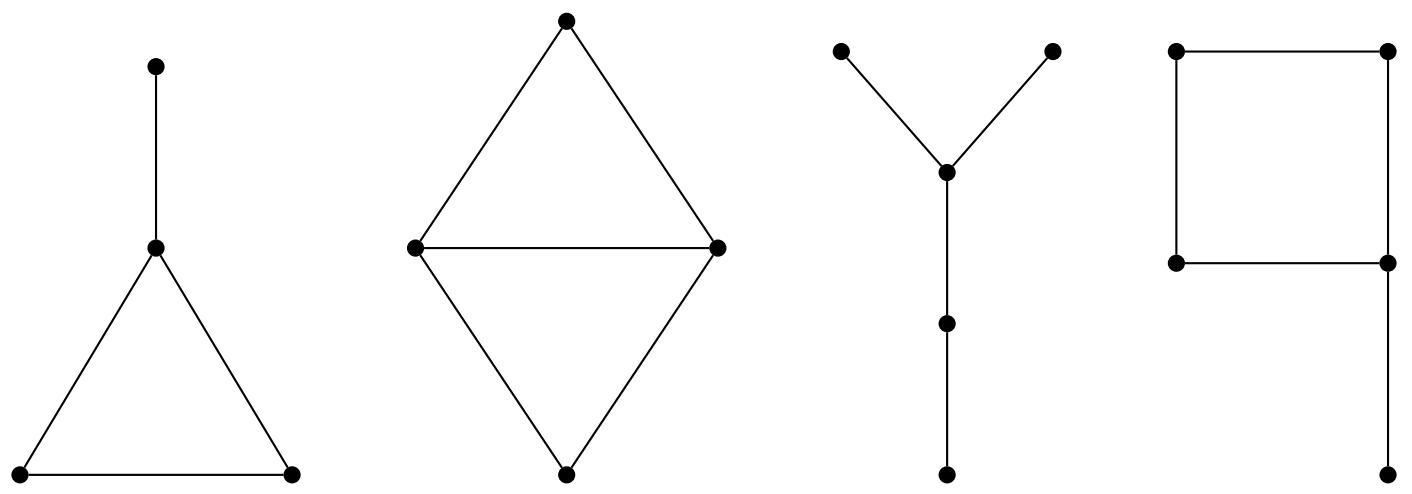

\begin{proposition}\label{prop:list-omega-1}
    Let $H$ be one of the following graphs: 
    \begin{enumerate}
        \item the paw graph;
        \item the diamond graph;
        \item the complete graph $K_4$;
        \item the 5-path graph $P_5$;
        \item the 5-cycle graph $C_5$;
        \item the fork graph;
        \item the banner graph.
    \end{enumerate} Then $\omega(H,\Bbbk)=1<2=\pdim I(H)$.
\end{proposition}

\begin{proof}
    One can compute $\HS_k(I(G))$ with Macaulay2 \cite{M2} for any integer $k$ and any $G$ of the seven graphs above, as all of them have at most five vertices. The computations in particular yield $\omega(H,\Bbbk)=1<2=\pdim I(H)$, as desired.
\end{proof}

We are now ready to prove the main theorem of this section. Recall that for any two integers $m,n\geq 1$, the \emph{$(m,n)$-complete bipartite graph}, denoted by $K_{m,n}$, is the graph with the vertex set $\{x_1,\dots, x_m, y_1, \dots, y_n\}$ and the edge set
\[
E(K_{m,n})=\{ \{x_i,y_j\} \mid  i
\in [m], j\in [n]\}.
\]
The graph $K_{1,n-1}$ is also called the \emph{$n$-star graph}.

\begin{theorem}\label{thm:omega-graphs}
    Let $G$ be a connected graph with at least one edge. Then one of the following holds:
    \begin{enumerate}
        \item $G$ has an induced $C_n$ for some $n\geq 3$ where $n\neq 4$ or an induced $P_4$.
        \item $G=K_{1.n}$ for some $n\geq 1$.
        \item $G=K_{m,n}$ for some $m,n\geq 2$.  
    \end{enumerate}
    We also have
    \[
    \omega(G,\Bbbk) = \begin{cases}
        1 & \text{if (1) holds},\\
        n-1 & \text{if (2) holds},\\
        2 & \text{if (3) holds}.
    \end{cases}
    \]
\end{theorem}

\begin{proof}
   Assume that $G$ is a tree, i.e., $G$ has no induced cycles. Then either $G$ is a star graph or $G$ has an induced $P_4$. Thus (1) or (2) holds, as desired.

    Now assume that $G$ has an induced cycle $C_n$ for some $n\geq 3$. Further, we assume that (1) does not hold. In particular, $G$ is a bipartite graph (a graph with no induced odd cycles) with no induced $P_4$. Statement (3) then follows from \cite[Lemma~2]{AlekseevLozin2004}.

    We now compute $\omega(G,\Bbbk)$. 
    
    If (2) holds, then the result follows from Lemma~\ref{lem:Taylor-omega}. 
    
    Now suppose that (1) holds, i.e., that $G$ has an induced cycle $C_n$ for some $n\neq 4$ or an induced $P_4$. If $n\geq 5$, then $G$ has an induced $C_5$ or $P_5$, and $\omega(G,\Bbbk)=1$ by Corollary~\ref{cor:omega=1}  and Proposition~\ref{prop:list-omega-1}, as desired. If $G$ has an induced $C_3$, since $G$ is connected, we have two cases:
    \begin{itemize}
        \item $G$ is isomorphic to $C_3$;
        \item $G$ has an induced paw graph, diamond graph, or $K_4$.
    \end{itemize}
    The former gives $\omega(G,\Bbbk)=1$ by direct computations with Macaulay2 \cite{M2}, and the latter gives $\omega(G,\Bbbk)=1$ from Corollary~\ref{cor:omega=1}  and Proposition~\ref{prop:list-omega-1}, as desired. Now we can assume that $G$ has an induced $P_4$. Similarly, since $G$ is connected, one of the following holds:
    \begin{itemize}
        \item $G$ is isomorphic to $P_4$;
        \item $G$ has an induced $P_5$, $C_5$, paw graph, diamond graph, fork graph, or banner graph.
    \end{itemize}
    By similar arguments as above, we also have $\omega(G,\Bbbk)=1$ in this case.

    Now we can assume that (3) holds, i.e., $G=K_{m,n}$ for some $m,n\geq 2$. To prove $\omega(G,\Bbbk)=2$, it suffices to show that
    \begin{equation}\label{eq:last-piece}
        \HS_2(I(G))\supseteq I(G)^{\{3\}}\quad \text{and} \quad \HS_3(I(G))\subsetneq I(G)^{\{4\}}.
    \end{equation}
    Set $V(G)=\{x_1,\dots, x_m, y_1, \dots, y_n\}$ and 
    \[
    E(K_{m,n})=\{ \{x_i,y_j\} \mid  i
    \in [m], j\in [n]\}.
    \]
    It is then straightforward that
    \begin{align*}
        I(G)^{\{3\}} = &(x_ix_jy_ky_l \mid 1\leq i<j \leq m \text{ and } 1\leq k<l\leq n) + \\
        &  (x_ix_jx_ky_l\mid 1\leq i<j <k\leq m \text{ and } l\in [n]) + \\
        &  ( y_iy_jy_kx_l\mid 1\leq i<j <k\leq n \text{ and } l\in [m]), \\
            I(G)^{\{4\}} =& (x_ix_jy_ky_l \mid 1\leq i<j \leq m \text{ and } 1\leq k<l\leq n) + \\
            &(x_ix_jx_kx_ly_r \mid 1\leq i<j<k<l \leq m \text{ and } r\in [n]) +\\
            &(y_iy_jy_ky_lx_r \mid 1\leq i<j<k<l \leq n \text{ and } r\in [m]).
    \end{align*}
    By Fr\"oberg's theorem \cite{Fro}, $I(G)$ has a linear resolution since the complement of $G$ is a disjoint union of two complete graphs, and thus chordal. In particular, $\HS_3(I(G))$ is generated in degree $5$. Combining this with the fact that $I(G)^{\{4\}}$ has a generator of degree $4$, we obtain
    \[
    \HS_3(I(G))\subsetneq I(G)^{\{4\}},
    \]
    which is the second part of (\ref{eq:last-piece}). We shall now show that the first part of (\ref{eq:last-piece}) also holds, concluding the proof.

    Fix integers $1\leq i<j\leq m$ and $1\leq k<l\leq n$. Let $H$ be the induced subgraph of $G$ with the vertex set $\{x_i,x_j,y_k,y_l\}$. The graph $H$ is a 4-cycle, and by Macaulay2 \cite{M2} and the restriction lemma (Lemma~\ref{lem:restriction-lemma}), we have
    \[
    x_ix_jy_ky_l \in \HS_2(I(H)) =\HS_2(I(G))^{\leq \prod_{x\in V(H)}x }\subseteq \HS_2(I(G)).
    \]
    By similar arguments, for any integers $1\leq i<j<k\leq m$ and $l\in [n]$, and $1\leq i'<j'<k'\leq n$ and $l'\in [m]$, we also have
    \[
    x_ix_jx_ky_l, y_{i'}y_{j'}y_{k'}x_{l'}\subseteq \HS_2(I(G)).
    \]
    Therefore $I(G)^{\{3\}}\subseteq \HS_2(I(G))$, as desired.
\end{proof}

\section{The CW-complex structure of the minimal resolution of $(x_1,\dots, x_q)^{\{r\}}$}

Given a CW-complex $\Delta$ with vertices labeled with monomials in a monomial set $M\subseteq S$, where $S$ is a polynomial ring over a field,  there is a natural way to induce a complex $\mathcal{F}_\Delta$ of free $S$-modules. We say that $\Delta$ \emph{supports a free resolution} of $(M)$ if $\mathcal{F}_{\Delta}$ is a free resolution of $(M)$, in which case $\mathcal{F}_\Delta$. For this section, the details of the construction of $\mathcal{F}_\Delta$ shall  not play a role, and instead we refer  to \cite[page 148]{BW02}. It is known that Morse resolutions are cellular (\cite[Theorem~1.3]{BW02}). This shall be our main source of cellular resolutions.

We remark that $(x_1,\dots, x_q)^r$ and $(x_1,\dots, x_{q+r-1})^{\{r\}}$ have the same total Betti numbers (\cite[Example~1]{EliKer} and Theorem~\ref{thm:extremal-Betti}). This raises a natural question.

\begin{question}\label{ques:same-CW-complex}
    Can the minimal free resolutions of the two ideals $(x_1,\dots, x_q)^r$ and $(x_1,\dots, x_{q+r-1})^{\{r\}}$ be supported on the same CW complex?
\end{question}

There are two problems with approaching this question. Firstly, the minimal resolution of an ideal can be supported by many CW-complexes. Secondly, there are virtually no tools for this type of problem, except for constructing such a CW-complex. We shall prove a sufficient~condition.

\begin{lemma}\label{lem:sufficient-condition-same-CW-complex}
    Let $I$ and $J$ be monomial ideals. Assume that the following holds:
    \begin{enumerate}
        \item $I$ and $J$ have the same total Betti numbers;
        \item there exists an lcm-preserving bijection  $2^{\mingens(I)}\to 2^{\mingens(J)}$;
        \item $I$ has a minimal Morse resolution.
    \end{enumerate}
    Then the minimal resolutions of $I$ and $J$ are supported on the same CW-complex.
\end{lemma}

We can now offer some positive partial answers to Question~\ref{ques:same-CW-complex}.

\begin{theorem}\label{thm:q=2-or-r=2}
    The minimal resolutions of the two ideals $(x_1,\dots, x_q)^r$ and $(x_1,\dots, x_{q+r-1})^{\{r\}}$ are supported on the same CW complex if either $q=2$ or $r=2$.
\end{theorem}

\begin{proof}
    It is known that the two ideals $(x_1,\dots, x_q)^r$ and $(x_1,\dots, x_{q+r-1})^{\{r\}}$ have the same total Betti numbers (\cite[Example 1]{EliKer}) and the former has a minimal Morse resolution (\cite[Section 5]{BW02}). By Lemma~\ref{lem:sufficient-condition-same-CW-complex}, it now suffices to construct an lcm-preserving bijection $2^{\mingens((x_1,\dots, x_q)^r)}\to 2^{\mingens((x_1,\dots, x_{q+r-1})^{\{r\}})}$, provided that $q=2$ or $r=2$.    

    We define two bijections.
\vskip.3truecm
The first is the map $\phi:2^{(x_1,\dots, x_n)^2}\to 2^{(y_1,\dots, y_{n+1})^{\lbrace 2\rbrace}}$ such that, for all indices $i,j\in\{1,\dots, n\}$,
\begin{eqnarray*}
\phi(x_i^2)&=&y_iy_{n+1}\\
\phi(x_ix_j)&=&y_iy_j.
\end{eqnarray*}

The second is the map $\psi:2^{(x_1,x_2)^n}\to 2^{(y_1,\dots, y_{n+1})^{\lbrace n\rbrace}}$ such that 
\begin{eqnarray*}
\psi(x_1^n)&=&y_1y_3\dots y_{n+1}\\
\psi(x_2^n)&=&y_2y_3\dots y_{n+1}\\
\psi(x_1^rx_2^{n-r})&=&y_1y_2y_3\dots y_{2r}\widehat{ y_{2r+1}}\dots y_{n+1}, \mbox{if $0<r< n-r$}, \\
\psi(x_1^{n-s}x_2^s)&=&y_1y_2y_3\dots y_{2s} y_{2s+1}\widehat{y_{2s+2}}\dots y_{n+1}, \mbox{if $0<s< n-s$},\\
\psi(x_1^rx_2^{n-r})&=&y_1y_2y_3\dots y_{n} , \mbox{if $r= n-r=\frac{n}2$}. \\
\end{eqnarray*}

    We shall show that $\phi$ is lcm-preserving. Consider $\sigma\in 2^{\mingens((x_1,\dots, x_n)^2)}$. It suffices to show that $\lcm(\phi(\sigma))$ is entirely determined by $\lcm(\sigma)$.  By relabelling if needed, we can set
    \[
    \lcm(\sigma)=x_{1}^2\cdots x_u^2x_{u+1}\cdots x_v
    \]
    for some $0\leq u\leq v\leq n$. We have two cases:
    \begin{itemize}
        \item \textbf{Case 1:} $u=0$, i.e., $\lcm(\sigma)=\prod_{i\in [v]}x_i$ is squarefree. This implies that $\sigma$ does not contain any square of a variable. By the definition of $\phi$, we have
        \[
        \lcm(\phi(\sigma))= \prod_{i\in [v]} y_i= \prod_{i\in [n], x_i\mid \lcm(\sigma)} y_i,
         \]
         dependent only on $\lcm(\sigma)$, as desired.
         \item \textbf{Case 2:} $u>0$. This guarantees that $y_{n+1}\mid \lcm(\phi(\sigma))$. By the definition of $\phi$, we have
        \[
        \lcm(\phi(\sigma))= y_{n+1} \prod_{i\in [v]} y_i= y_{n+1}\prod_{i\in [n], x_i\mid \lcm(\sigma)} y_i,
         \]
         dependent only on $\lcm(\sigma)$, as desired.
    \end{itemize}
    It remains to show that $\psi$ is lcm-preserving. Consider $\sigma\in 2^{\mingens((x_1, x_2)^n)}$. It suffices to show that $\lcm(\psi(\sigma))$ is entirely determined by $\lcm(\sigma)$. Indeed, if $\lcm(\sigma)\in \mingens((x_1, x_2)^n)$, then $|\sigma|=1$ and thus $\lcm(\psi(\sigma))$ is uniquely determined by $\lcm(\sigma)$. On the other hand, if $\lcm(\sigma)\notin \mingens((x_1, x_2)^n)$, i.e., $|\sigma|\geq 2$, then $\lcm(\psi(\sigma))=y_1\cdots y_{n+1}$. This concludes the proof.
\end{proof}

\begin{theorem}\label{thm:impossible-lcm-map}
    Between $\mingens((x_1,x_2, x_3)^3) $ and $\mingens ((y_1,y_2, y_3,y_4,y_5)^{\{3\}})$, no lcm-preserving map exists.
\end{theorem}

\begin{proof}
    Set
    \[
    X=\mingens((x_1,x_2, x_3)^3)= \{x_1^3,x_1^2x_2,x_1^2x_3,x_1x_2^2,x_1x_2x_3,x_1x_3^2,
x_2^3,x_2^2x_3,x_2x_3^2,x_3^3\}
    \]
    and
\[
Y=\mingens ((y_1,y_2, y_3,y_4,y_5)^{\{3\}})=\{y_i y_j y_k : 1\le i<j<k\le 5\}.
\]
We note that there is no bijection from $Y$ to $X$ such that the power map is lcm-preserving, as the LCM lattice of $(y_1,y_2, y_3,y_4,y_5)^{\{3\}}$ has less elements than that of  $(x_1,x_2, x_3)^3$.
Suppose by contradiction that there is a bijection $\psi:X\to Y$ such that the power map  is lcm-preserving. In particular, for all $m_1,m_2, m'_1,m'_2\in X$, whenever
\[
\mathrm{lcm}(m_1,m_2)=\mathrm{lcm}(m'_1,m'_2),
\]
we also have
\[
\mathrm{lcm}(\psi(m_1),\psi(m_2))
=
\mathrm{lcm}(\psi(m'_1),\psi(m'_2)).
\]
Now consider the complete graph $K_5$ on the vertex set $[5]=\lbrace1,\dots,5\}$ and define the bijection $F:Y\to E(K_5)$ by setting  $F(y_iy_jy_k)=[5]\setminus\lbrace i,j,k\rbrace$ for all triples $(i,j,k)$. 

Then consider the following three subsets of $X$:
\[
X_1=\{x_2x_3^2,\ x_1x_3^2,\ x_1x_2x_3\},
\]
\[
X_2=\{x_2^2x_3,\ x_1x_2^2,\ x_1x_2x_3\},
\]
\[
X_3=\{x_1^2x_3,\ x_1^2x_2,\ x_1x_2x_3\}.
\]

Fix an index $1\leq i\leq 3$. In $X_i$, all pairs of elements have the same lcm. Therefore, the elements of $X_i$ are mapped by $F\circ \psi$ to three edges having the same pairwise intersection, i.e., having a common vertex. In other words, $F\circ \psi(X_i)$ is a $3-$edge star graph $S_i$.

Since  $X_1, X_2, X_3$ have the monomial \(x_1x_2x_3\) in common, $S_1, S_2, S_3$ must have (exactly) one edge $e$ in common. But then the center of each $S_i$ is an endpoint of $e$. Hence at most two different centers are possible. Thus two stars, say $S_1$ and $S_2$, must share the same center, say $1$. Now, $S_1$ and $S_2$ have two edges besides the common edge $e$, say $\lbrace 1,2\rbrace$, $\lbrace 1,3\rbrace$ and  $\lbrace 1,4\rbrace$, $\lbrace 1,5\rbrace$, respectively. Since also the edge $e$ is incident to the vertex $1$, this provides a contradiction.

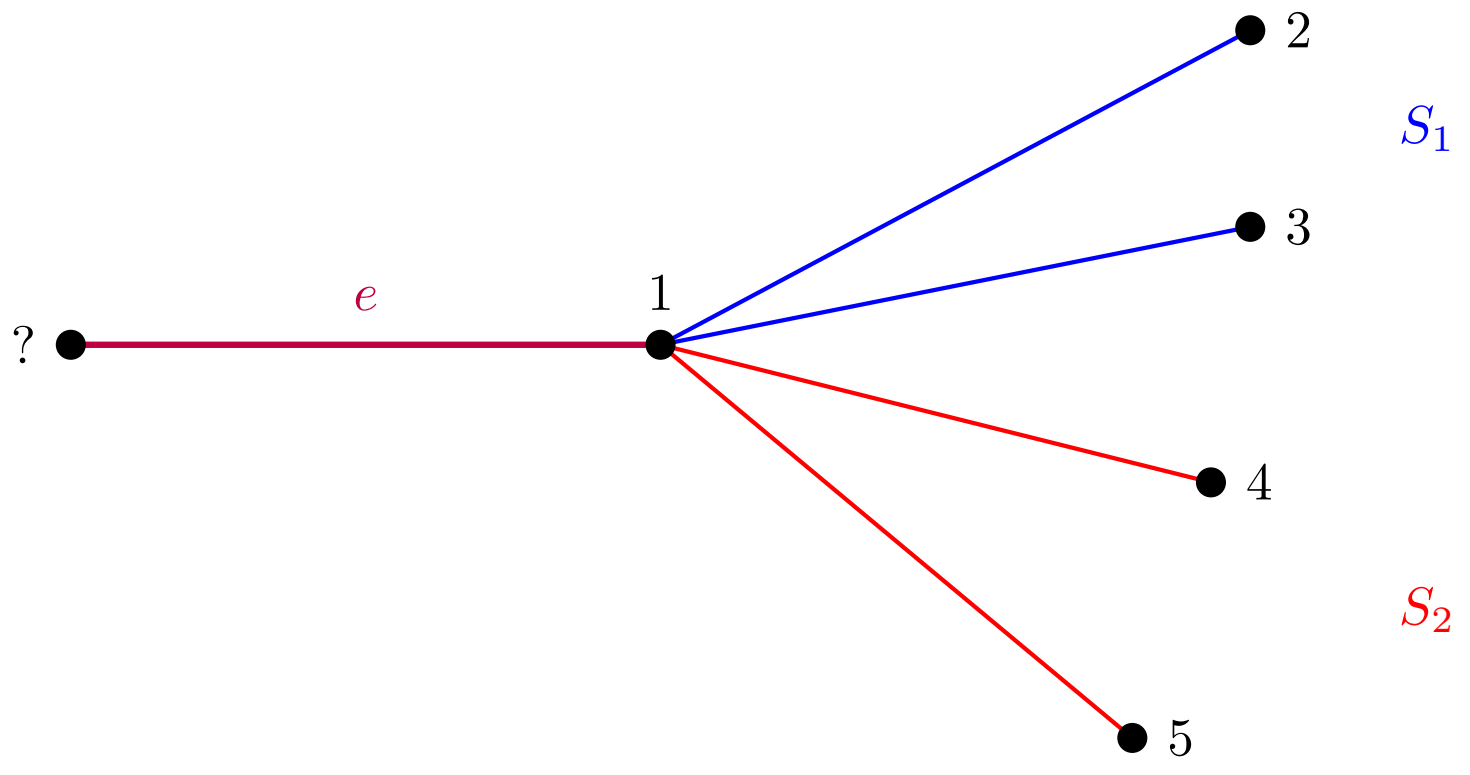
\vspace{1em}

\begin{center}
\begin{tikzpicture}[scale=1.3,
    vertex/.style={circle, fill=black, inner sep=2pt}
]

\coordinate (v1) at (0,0);
\coordinate (v2) at (3,1.6);
\coordinate (v3) at (3,0.6);
\coordinate (v4) at (2.8,-0.7);
\coordinate (v5) at (2.4,-2);
\coordinate (v6) at (-3,0);

\draw[blue, thick] (v1)--(v2);
\draw[blue, thick] (v1)--(v3);

\draw[red, thick] (v1)--(v4);
\draw[red, thick] (v1)--(v5);

\draw[purple, very thick] (v1)--(v6);
\node[purple, above=3pt] at (-1.5,0) {$e$};
\node[vertex,label=above:$1$] at (v1) {};
\node[vertex,label=right:$2$] at (v2) {};
\node[vertex,label=right:$3$] at (v3) {};
\node[vertex,label=right:$4$] at (v4) {};
\node[vertex,label=right:$5$] at (v5) {};
\node[vertex,label=left:$?$] at (v6) {};

\node[blue] at (3.9,1.1) {$S_1$};
\node[red]  at (3.9,-1.35) {$S_2$};

\end{tikzpicture}
\end{center}

\vspace{1em}

\noindent
Therefore, there is no bijection \(\psi\) satisfying the required condition.
\end{proof}

\begin{remark}
    Theorem~\ref{thm:impossible-lcm-map} does not provide a negative answer to Question~\ref{ques:same-CW-complex} in the case $q=r=3$ (and also in the cases where $q\geq 3$ and $r\geq 3$). It simply means that our sufficent condition (Lemma~\ref{lem:sufficient-condition-same-CW-complex}) is rendered ineffective for higher cases, which we shall leave open.
\end{remark}


\bibliographystyle{amsplain}
\bibliography{refs}
\end{document}